\documentclass[12pt]{article}
\usepackage[reqno]{amsmath}
\usepackage{amssymb, theorem, enumerate}
\usepackage{graphicx}

\expandafter\ifx\csname pdfpagebox\endcsname\relax\else
\fi

\theoremstyle{change}
{\theorembodyfont{\slshape}
\newtheorem{theorem}{Theorem.}[section]
\newtheorem{lemma}[theorem]{Lemma.}
\newtheorem{corollary}[theorem]{Corollary.}}
\theorembodyfont{\rmfamily}

\newcommand\lref[1]{Lemma~\ref{lem:#1}}
\newcommand\tref[1]{Theorem~\ref{thm:#1}}
\newcommand\cref[1]{Corollary~\ref{cor:#1}}
\newcommand\sref[1]{Section~\ref{sec:#1}}

\def\proof{\noindent{{\sl Proof. }}}
\def\sqr#1#2{{\vbox{\hrule height.#2pt
    \hbox{\vrule width.#2pt height#1pt \kern#1pt
        \vrule width.#2pt}\hrule height.#2pt}}}
\def\eqed{\sqr53}
\def\qed{%
    \ifmmode\eqno\eqed
    \else\nobreak\ \hfill\eqed\medbreak\fi}

\newcommand\al{\alpha}
\newcommand\be{\beta}

\newcommand\la{\lambda}

\renewcommand\th{\theta} %--- Latex uses \th for a Norse character

\newcommand\fld{{\mathbb F}}

\newcommand\re{{\mathbb R}}%reals

\newcommand\comp[1]{{\mkern2mu\overline{\mkern-2mu#1}}}
\newcommand\diff{\mathbin{\mkern-1.5mu\setminus\mkern-1.5mu}}% for \setminus
\newcommand\sbs{\subseteq}

\newcommand\seq[3]{#1_{#2},\ldots,#1_{#3}}
\DeclareMathOperator{\supp}{supp}

\newcommand\pmat[1]{\begin{pmatrix} #1 \end{pmatrix}}

\newcommand\one{{\bf1}}
\DeclareMathOperator{\rk}{rk}

\DeclareMathOperator{\col}{col}
\newcommand\mat[3]{\mathrm{Mat}_{#1\times #2}(#3)}

\DeclareMathOperator\dist{dist}

\newcommand\what{\widehat}

\title{Controllable Subsets in Graphs}
\author{Chris Godsil}

\begin{document}
\maketitle

\begin{abstract}
	Let $X$ be a graph on $v$ vertices with adjacency matrix $A$, and let
	let $S$ be a subset of its vertices with characteristic vector $z$.
	We say that the pair $(X,S)$ is \textsl{controllable} if the vectors
	$A^rz$ for $r=1,\ldots,v-1$ span $\re^v$.  Our concern is chiefly with the
	cases where $S=V(X)$, or $S$ is a single vertex.  In this paper we develop 
	the basic theory of controllable pairs.  We will see that if $(X,S)$ is 
	controllable then the only automorphism of $X$ that fixes $S$ as a set 
	is the identity.  If $(X,S)$ is controllable for some subset $S$ then 
	the eigenvalues of $A$ are all simple.
\end{abstract}

\section{Introduction}

Let $X$ be a graph on $v$ vertices with adjacency matrix $A$.  If $z\in\re^v$,
define the matrix $W$ by
\[
	W= \pmat{z&Az&\hdots&A^{v-1}z}
\]
The pair $(A,z)$ is \textsl{controllable} if $W$ is invertible.  In this article,
$z$ will often be the characteristic vector of some subset $S$ of $V(X)$, and
then we will say that $(X,S)$ is controllable if $(A,z)$ is.  When $S\sbs V(X)$, 
the entries of $W=W_S$ counts walks in the graph $X$, and we will call $W_S$ 
the \textsl{walk matrix} of $S$.

We note one interesting property of controllable pairs.

\begin{lemma}\label{lem:aut}
	If $(X,S)$ is controllable, then any automorphism of $X$ that fixes $S$ as a
	set is the identity.
\end{lemma}

\proof
We view the automorphisms of $X$ as permutation matrices that commute with $A$.
Let $z$ be the characteristic vector of $S$.
An automorphism of $X$ fixes $S$ if and only if $Pz=z$.  If $Pz=z$ then 
for $r=0,\ldots,v-1$
\[
	PA^r z =A^r Pz =A^r z
\] 
and therefore $PW_S=W_S$.  Hence if $W_S$ is invertible, $P=I$.\qed

In this paper we develop the theory of controllable pairs. We will see that there
is a close connection to the subject of control theory. The ideas in this paper have
already been put to use in quantum physics---see \cite{MeSev}.

\section{Characterizations of Controllability}
\label{sec:ChaCon}

We derive some useful characterizations of controllability.

Assume that $A=A(X)$ has the spectral decomposition
\[
	A =\sum_\th \th E_\th.
\]
Then
\[
	A^r z = \sum_\th \th^r E_\th z
\]
and hence $\col(W)$ is spanned by the vectors $E_\th z$, that is, by the nonzero
projections of $z$ onto the distinct eigenspaces of $A$. We say that an
eigenvalue $\th$ is in the \textsl{support} of $y$ if $E_\th y\ne0$.
Equivalently, $\th$ is in the support of $y$ if $y^TE_\th y\ne0$. The
\textsl{dual degree} of $y$ is $|\supp(y)| -1$. So $(X,S)$ is controllable if and
only if the dual degree of $S$ is $v-1$. Note that if $X$ is connected and $z$ is
the eigenvector of $X$ with eigenvalue equal to the spectral radius of $A$, then
by the Perrron-Frobenius theory, all entries of $z$ are positive. It follows that
the dual degree is non-negative.

From the spectral decomposition of $A$, we see that
\[
	(tI-A)^{-1}z = \sum_\th \frac1{t-\th} E_\th z
\]
and hence
\[
	z^T(tI-A)^{-1}z = \sum_\th \frac{z^TE_\th z}{t-\th}
\]  
Since 
\begin{equation}
	\label{eq:ratfun}
	z^T E_\th z = z^T E_\th^2 z =(E\th z)^T E_\th x
\end{equation}
we have that $z^TE_\th z=0$ if and only if $E_\th z=0$.  Therefore the rank of $W_S$ is
equal to the number of poles of the rational function $z^T(tI-A)^{-1}z$.
There is a polynomial $\phi_S(X,t)$ with degree at most $v-1$ such that
\[
	z^T(tI-A)^{-1}z =\frac{\phi_S(X,t)}{\phi(X,t)}
\]
(It is not hard to show that, if $S$ is the vertex $u$, then 
$\phi_S(X,t)=\phi(X\diff u,t)$.)

This provides a useful characterization of controllability:

\begin{lemma}\label{lem:coprime}
	Let $X$ be a graph on $n$ vertices and suppose $S\sbs V(X)$.  Let $z$ be the
	characteristic vector of $S$.  Then $(X,S)$ is controllable if and only if 
	the rational function $z^T(tI-A)^{-1}z$ has $n$ distinct poles.\qed
\end{lemma}

Our next result characterizes controllability in terms of linear algebra rather
than rational functions,

\begin{theorem}\label{thm:algebra}
	Let $S$ be a subset of the vertices of the graph $X$, with
	characteristic vector $z$. The following statements are equivalent:
	\begin{enumerate}[(a)]
		\item 
		$(X,S)$ is controllable.
		\item
		The matrices $A$ and $zz^T$ generate the algebra of all 
		$v\times v$ matrices.
		\item
		The matrices $A^i zz^T A^j$ where $0\le i,j < v$ form a basis for the
		algebra of all $v\times v$ matrices.
	\end{enumerate}
\end{theorem}

\proof
We show that (a) and (c) are equivalent.

If $(X,S)$ is controllable, then the vectors
\[
	z, Az,\ldots, A^{v-1}z
\]
are linearly independent in $\re^v$.  Since
\[
	\re^v\otimes\re^v \cong \mat vv\re
\]
it follows that the matrices 
\[
	(A^i z)(A^j z)^T = A^i zz^T A^j,\qquad(0\le i,j< v)
\]
are linearly independent in $\mat vv\re$.  On the other hand, if $(X,S)$ is not
controllable, then the vectors $A^rz$ span a space of dimension at most $v-1$,
and the matrices $A^i zz^T A^j$ span a space of dimension at most $(v-1)^2$.

To complete the proof, note that
\[
	zz^T A^r zz^T =(z^TA^rz)zz^T
\]
and therefore any element of the algebra generated by $A$ and $zz^T$ is a linear
combination of matrices of the form
\[
	A^r zz^T A^s,\qquad(r,s\ge0)
\]
Since $A$ is $v\times v$, any polynomial in $A$ is a linear combination of the
powers
\[
	I,A,\ldots,A^{v-1},
\]
We conclude that (b) implies (c).  Since (b) is an immediate consequence of (c),
we are done.\qed

\section{Isomorphism}

In this section we consider pairs that need not be controllable.

Let $X$ be a graph on $v$ vertices with adjacency matrix $A$ and let $y$ be a vector
in $\re^v$.  Let $Y$ be a graph on $v$ vertices with adjacency matrix $B$ and
let $z$ be a vector in $\re^v$. We say that the pairs $(X,y)$ and $(Y,z)$
are \textsl{isomorphic} if there is an orthogonal matrix $L$ such that
\[
	LAL^T =B,\quad Ly=z.
\]
In this case $LW_y=W_z$; thus controllability is preserved by isomorphism.
Further 
\[
	W_z^TW_z =W_y^TL^TLW_y =W_y^TW_y.
\]
We will occasionally refer to the characteristic polynomial of the matrix $A$ as
the characteristic polynomial of the pair $(A,z)$. 

If the pairs $(A,y)$ and $(B,z)$ are isomorphic, then $A$ and $B$ must have the same
characteristic polynomial.  

\begin{theorem}
	\label{thm:ratfun}
	Two pairs $(A,y)$ and $(B,z)$ are isomorphic if and only if $A$ and $B$
	are similar and
	\[
		y^T(I-tA)^{-1}y =z^T(I-tB)^{-1}z.
	\]
\end{theorem}

\proof
We have seen that the necessity of this condition is an easy consequence
of the definition.  So we assume that $A$ and $B$ are cospectral and
that our two rational functions are equal.  From our remarks at the start
of \sref{ChaCon}, in particular \eqref{eq:ratfun}, the latter condition implies that 
$y$ and $z$ have the same support
and that $y^TE_\th y=z^TF_\th z$ for each eigenvalue $\th$ in $\supp(y)$.  

We construct two orthonormal bases for $\re^v$;
the linear map that takes the first basis to to the second will be our isomorphism.

Let
\[
	A=\sum_\th \th E_\th,\qquad B=\sum_\th \th F_\th
\]
be the spectral decompositions of $A$ and $B$.  

Construct an orthonormal basis $\al$ for
$\re^v$ as follows.  Suppose $\rk(W_y)=d$.  The first $d$ vectors of the
basis will be the normalizations of the non-zero vectors $E_\th y$.
For each $\th$ in $\supp(y)$, add an orthonormal basis for the
subspace of eigenvectors for $A$ with eigenvalue $\th$ that are orthogonal to $E_\th y$.
If $\th\notin\supp(y)$, add an orthonormal basis for the $\th$-eigenspace of $A$.
By the same procedure we can form an orthonormal basis $\be$ relative to $B$ and $y$
and, possibly after some rearrangement, we may assume that the vectors $\al_i$
and $\be_i$ have the same eigenvalue for each $i$.  If $L$ is the matrix representing
the unique linear mapping that sends $\al$ to $\be$, then $L$ is orthogonal
and $LAL^T=B$.

Set $W=W_y$ and let $Y$ be the matrix whose columns are the nonzero vectors $E_\th y$.
Then $W$ and $Y$ have the same column space, in fact if $M$ is the $d\times v$
matrix whose $ij$-entry is $\th^{j-1}$, where $\th$ is the $i$-th eigenvalue 
in $\supp(y)$, then $W=YV$.  Note that $V$ is determined by the support of $y$.
The columns of $Y$ are pairwise orthogonal and therefore
\[
	W^TW = V^T DV,
\]
where $D$ is the diagonal matrix with diagonal entries of the form $y^TE_\th y$.
Since $y^T(I-tA)^{-1}y =z^T(I-tB)^{-1}z$, we infer that $W_y^TW_y=W_z^TW)_z$
and, since $V$ has a right inverse, we conclude that
\[
	y^TE_\th y = z^TF_\th z.
\]
for each eigenvalue $\th$ in $\supp(y)$.  Since
\[
	y =\sum_{\th\in\supp(y)} E_\th y
\]
and since $L$ maps $(y^TE_\th y)^{-1/2}E_\th y$ to $(z^TF_\th z)^{-1/2}F_\th z$,
it follows that $Ly=z$.\qed

Note that the two rational functions above are equal if and only if
\[
	y^T E_\th y = z^T F_\th z
\]
for \textbf{all} eigenvalues $\th$.

\begin{corollary}
	If $(A,y)$ and $(B,y)$ are controllable and $y^T(I-tA)^{-1}y=z^T(I-tB)^{-1}z$,
	then $(A,y)$ and $(B,z)$ are isomorphic.
\end{corollary}

\proof
If $(A,y)$ is controllable, then the eigenvalues of $A$ are distinct and each one
is s pole of $y^T(I-tA)^{-1}y$.  So our hypothesis implies that $A$ and $B$
are cospectral and that $y^T E_\th y=z^TF_\th z$ for all eigenvalues $\th$.\qed

By Lemma~2.2 in \cite{GMrecon} (for example) it follows that if $X$ and $Y$
are cospectral then $\comp{X}$ and $\comp{Y}$ are cospectral if and only if
\[
    \one^T(I-tA(X))^{-1}\one = \one^T(I-tA(Y))^{-1}\one.
\]
So the results of this section imply the important result of Johnson
and Newman \cite{JohnNew} that if $X$ and $Y$ are cospectral with cospectral
complements, then there is an orthogonal matrix $L$ such that
\[
    L^TA(X)L=A(Y),\quad L^T(A(\comp{X}))L = A(\comp{Y}).
\]

\section{Graph Theory}

If $S\sbs V(X)$, we define the \textsl{covering radius} of $S$ to be the least 
integer $r$ such that each vertex of $X$ is at distance at most $r$ from a vertex
of $S$.  Thus $S$ has covering radius equal to 1 if and only if it is a dominating set,
and the diameter of $X$ is the maximum value of the covering radius of a vertex.

\begin{lemma}
	\label{lem:covrad}
	If $S$ has dual degree $m$ and covering radius $r$, then $r\le m$.
\end{lemma}

\proof
If $v\in V(X)$, then $(A^iz)_v$ is equal to the number of walks of length $i$
from $v$ to a vertex in $S$.  It follows that $((A+I)^iz)_v$ is zero if and only
if $i$ is less than $\dist(v,S)$.  From this it follows in turn that the vectors
\[
	z,Az,\ldots,A^rz
\]
are linearly independent and therefore $r+1$ is a lower bound on $\rk(W_z)$.\qed

One consequence of this lemma is the well known result that if $X$ has diameter $d$, 
then $d+1$ is less than or equal to the number of distinct eigenvalues of $A$.
As an example, if $X$ is the path $P_n$ on $n$ vertices and $S$ is one of its
end-vertices, then covering radius of $S$ is $n-1$.  Hence the dual degree
of $S$ is $n-1$, from which we deduce the well known fact that the eigenvalues 
of the path are distinct.

\begin{lemma}
	If $X$ is vertex transitive and $|V(X)|>2$, no subset of $V(X)$ is controllable.
\end{lemma}

\proof
If $X$ has a controllable subset with characteristic vector $z$, then
the vectors $E_\th z$ form a basic for $\re^v$, and thus $X$ has $v$ simple
eigenvalues.  But the only vertex transitive graph with all eigenvalues simple 
is $K_2$.\qed

If $S\sbs V(X)$, we define the \textsl{cone of X relative to $S$} to be the graph
we get by taking one new vertex and declaring it to be adjacent to each vertex in $S$.

\begin{theorem}
	\label{thm:cone}
	The pairs $(X,S)$ and $(Y,T)$ are isomorphic if and only if $X$ is cospectral
	to $Y$ and the cone of $X$ relative to $S$ is cospectral to the cone of $Y$
	relative to $T$.
\end{theorem}

\proof
Let $b$ denote the characteristic vector of $S$ and let $\what X$ denote the cone 
over $X$ relative to $S$.  Then
\[
	A(\what X)=\pmat{0&b^T\\ b&A}
\]
and so
\[
	\pmat{t& -b^T\\-b&tI-A} 	
		=\pmat{1&0\\0&tI-A} \pmat{t&-b^T\\ -(tI-A)^{-1}b&I}.
\]
Accordingly
\[
	\phi(\what X,t) =\phi(X,t)(t-b^T(tI-A)^{-1}b)
\]
and this yields that
\[
	\frac{\phi(\what X,t)}{\phi(X,t)} = t -\sum_\th \frac{b^TE_\th b}{t-\th}.
\]
Our result follows now from \tref{ratfun}.\qed

\begin{theorem}
	Suppose $V(X)=\{1,\ldots,v\}$ and $S\sbs V(X)$.  Construct the cone
	$\what X$ by joining the vertex $0$ to each vertex in $S$.  Then $(X,S)$
	is controllable if and only if $(\what{X},\{0\})$ is controllable.
\end{theorem}

\proof
Assume $v=|V(X)|$.
If $b$ is the characteristic vector of $S$, we have
\begin{equation}
	\label{eq:phwhat}
	\frac{\phi(\what X,t)}{\phi(X,t)} = t -\sum_\th \frac{b^TE_\th b}{t-\th}.
\end{equation}
Further $(X,S)$ is controllable if and only if this rational function has
$v$ distinct poles.

Now
\[
	e_0^T(tI-\what{A})^{-1}e_0 =\bigl((tI-\what{A})^{-1}\bigr)_{0,0}
		= \frac{\phi(X,t)}{\phi(\what{X},t)}
\]
and therefore $(\what{X},\{0\})$ is controllable if and only if the rational
function $\phi(X,t)/\phi(\what{X},t)$ has $v+1$ distinct poles, that is, if and
only if $\phi(\what{X},t)/\phi(X,t)$ has exactly $v+1$ distinct zeros.
Since the derivative of the right side in \eqref{eq:phwhat} is positive everywhere it is
defined, between each pair of consecutive zeros there is exactly one pole.  Therefore
there are $v+1$ distinct zeros.

The following corollary provides infinite families of controllable pairs.

\begin{corollary}
	Let $S$ be a subset of $V(X)$, and let $Y_k$ be the graph obtained by 
	taking a path on $k$ vertices and joining one of its end-vertices to 
	each vertex in $S$. Let $0$ denote the other end-vertex of the path.  
	If $(X,S)$ is controllable then $(Y_k,\{0\})$ is controllable.\qed
\end{corollary}

Our next result generalizes Lemma~2.4 from \cite{WangXu}.

\begin{lemma}
\label{lem:wtwsi}
	Suppose the pairs $(X,S)$ and $(Y,T)$ are isomorphic and controllable.
	Then the matrix $W_TW_S^{-1}$ represents the isomorphism from $(X,S)$
	to $(Y,T)$.
\end{lemma}

\proof
Let $A$ and $B$ be the adjacency matrices of $X$ and $Y$ respectively.

Since the pairs are isomorphic, $W_S^TW_S=W_T^TW_T$.  Since they are controllable,
$W_S$ and $W_T$ are invertible and therefore
\[
	W_TW_S^{-1} =W_T^{-T}W_S^T =(W_TW_S^{-1})^{-T} 
\]
Hence $Q+W_TW_S^{-1}$ is orthogonal.

Let $C$ denote the companion matrix of $\phi(X,t)$.  Then
\[
	AW_S = W_SC
\]
and, since $A$ and $B$ are similar,
\[
	BW_T = W_TC.
\]
Hence
\[
	BW_TW_S^{-1} = W_TCW_S^{-1} =W_TW_S^{-1} A
\]
and thus $B= QAQ^{-1}$.

Let $y$ and $z$ be the characteristic vectors of $S$ and $T$ respectively.
Since $QW_S =W_T$, we certainly have $Qy=z$.\qed

\begin{corollary}
	If the pairs $(X,S)$ and $(X,T)$ are isomorphic and controllable
	and $Q=W_TW_S^{-1}$, then $Q$ commutes with $A(X)$ and $Q^2=I.$
\end{corollary}

\proof
From the lemma we have $QAQ^{-1}=A$, so $Q$ and $A$ commute.  Since the eigenvalues
of $A$ are all simple, this implies that $Q$ is a polynomial in $A$ and therefore it
is symmetric.\qed

When the hypotheses of this corollary hold, the matrix $Q$ can be viewed as a
kind of ``approximate'' automorphism of order two---it is rational, commutes with
$A$ and swaps the characteristic vectors of $S$ and $T$. If $S$ and $T$ are single
vertices $u$ and $v$, then $Q$ will be block diagonal with one block of the form
\[
	\pmat{0&1\\ 1&0}
\]
and the other an orthogonal matrix of order $(v-2)\times(v-2)$ which commutes with the
adjacency matrix of $X\diff\{u,v\}$.

\section{Controllable Graphs}

We say that graph is \textsl{controllable} if $(X,V(X))$ is controllable. Since
any automorphism of $X$ fixes $V(X)$, we see that a controllable graph is
asymmetric. We can see this another way. If $(X,V(X))$ is controllable,
then $W$ is invertible and so if $e_u^TW =e_v^TW$, then $u=v$.
One consequence of this observation is that the ordering of the vertices
obtained from the lexicographic ordering of the rows of $W$ is canonical: two
controllable graphs are isomorphic if and only their ordered walk matrices 
are equal.

It is also immediate that a graph is controllable if and only if its complement
is.

\begin{theorem}
	\label{thm:2dcontrol}
	If $u$ and $v$ are cospectral vertices in $X$, then the $A$-modules
	generated by $e_u+e_v$ and $e_u-e_v$ are orthogonal.  If the $A$-module 
	generated by $\{e_u,e_v\}$ is $\re^{V(X)}$, then $\re^{V(X)}$ is the 
	direct sum of these two cyclic modules.
\end{theorem}

\proof
If $u$ and $v$ are cospectral, then $(E_\th)_{u,u}=(E_\th)_{v,v}$ for each
eigenvalue $\th$ of $X$. For any projection $E_\th$ we have
\[
	(e_u+e_v)^T E_\th(e_u-e_v) = (E_\th)_{u,u} -(E_\th)_{v,v}
\]
and so the vectors $E_\th(e_u+e_v)$ are orthogonal to the vectors
$E_\tau(e_u-e_v)$, for all choices of $\th$ and $\tau$.\qed

The second condition in the theorem will hold if $u$ (or $v$) is controllable.
The theorem implies that, if $u$ and $v$ are cospectral and $z$ lies in the
$A$-module generated by $e_u+e_v$, then $z_u=z_v$.

We have the following consequence of \lref{coprime} and the remark preceding it:

\begin{lemma}\label{lem:}
	A vertex $u$ in $X$ is controllable if and only if $\phi(X\diff u,t)$ 
	and $\phi(X,t)$ are coprime.\qed
\end{lemma}

For the path $P_n$ on $n$ vertices we have
\[
	\phi(P_0,t)=1,\qquad \phi(P_1,t)=t
\]
and, if $n\ge1$,
\[
	\phi(P_{n+1},t) =t\phi(P_n,t) -\phi(P_{n-1},t)
\]
from which it follows by induction that $\phi(P_{n+1},t)$ and $\phi(P_n,t)$ are
coprime for all $n$. So if $1$ is an end-vertex of $P_n$, the pair $(P_n,\{1\})$
is controllable.

\begin{corollary}\label{cor:}
    If the characteristic polynomial of $X$ is irreducible over the rationals,
    then $(X,V(X))$ is controllable and $(X,u)$ is controllable for any
    vertex $u$.\qed
\end{corollary}

In \cite{GMrecon} it is proved that controllable graphs are reconstructible.
We conjecture that almost all graphs are controllable.

\section{Laplacians}

The theory we have presented will hold for any symmetric matrix.  If $D$ is 
the diagonal matrix of valencies of the vertices of $X$, then
\[
	L := D-A
\]
is the \textsl{Laplacian} of $X$.  This is a symmetric matrix with row sums zero.
If $e_i$ and $e_j$ are two of the standard basis vectors, then
\[
	H_{i,j} := (e_i-e_j)(e_i-e_j)^T
\]
If the graph $Y$ is obtained by adding the edge $ij$ to $X$, then
\[
	L(Y) =D-A+H_{i,j}.
\]
Thus 
\[
	L(X) =\sum_{ij\in E(X)} H_{i,j}.
\]

Now
\begin{align*}
	\det(tI -L -H_{i,j}) &=\det[(tI-L)(I -(tI-L)^{-1})H_{i,j}]\\
		&= \det(tI-L) \det(I -(tI-L)^{-1}(e_i-e_j)(e_i-e_j)^T)\\
		&= \det(tI-L) (1-(e_i-e_j)^T(tI-L)^{-1}(e_i-e_j))
\end{align*}
and if $h:=e_i-e_j$, then
\[
	\frac{\phi(L(Y),t)}{\phi(L(X),t)}
		=1-h^T(tI-L)^{-1}h = 1 -\sum_\la \frac{h^T F_\la h}{t-\la}
\]
where $L=\sum_\la \la F_\la$ is the spectral decomposition of $L$. It follows
that the eigenvalues of $L(Y)$ are determined by the eigenvalues of $L(X)$ along
with the squared lengths of the projections of $e_i-e_j$ onto the eigenspaces of
$L(X)$.

If we get $Y$ from $X$ by deleting the edge $ij$, then we find that
\[
	\frac{\phi(L(Y),t)}{\phi(L(X),t)} = 1 +\sum_\la \frac{h^T F_\la h}{t-\la}
\]
Let $h$ denote $e_i-e_j$.  

We observe that $A^rh$ is orthogonal to $\one$, and so the dimension of the
$A$-module generated by $h$ is at most $v-1$. We say that the pair of vertices
$\{i,j\}$ is \textsl{controllable} relative to the Laplacian if
\[
	\pmat{\one&h&Lh&\hdots&L^{v-1}h}
\]
has rank $v-1$.  

If $ij$ is controllable and $P$ is an automorphism of $X$ that fixes 
$\{i,j\}$, then either 
\[
	P(e_i-e_j) =e_i-e_j
\]
and $PW=W$, or
\[
	P(e_i-e_j) =e_j-e_i
\]
and $PW=-W$.  In the latter case $P=-I$ and so it is not a permutation matrix, in
the former case $P=I$.  We conclude that if $ij$ is controllable, then only
the identity automorphism fixes the set $ij$.

\section{Control Theory}

In this section we provide a brief introduction to some concepts from
control theory. Our favorite source for this material is the book of 
Kailath \cite{Kailath} (but there is a lot of choice).

Consider a discrete system whose state at time $n$ is $x_n$, where
$x_n\in\fld^d$. The states are related by the recurrence
\begin{equation}
	\label{eq:control}
	x_{n+1} =Ax_n +u_n B \qquad(n\ge0).
\end{equation}
where $A$ and $b$ are fixed matrices and the $(u_n)_{n\ge0}$ is arbitrary.  
The output $c_n$ at time $n$ is equal to $c^Tx_n$, where
$c$ is fixed.  The basic problem is determine information about the state of the system
given $(u_n)$ and $(c_n)$.  From \eqref{eq:control} we find that
\[
	\sum_{n\ge0} t^nx_{n+1} =A\sum_{n\ge0}t^nx_n +\Bigl(\sum_{n\ge0}u_n t^n\Bigr)b
\]
If we define
\[
	X(t) :=\sum_{n\ge0}t^nx_n,\quad u(t) := \sum_{n\ge0}u_n t^n, 
		\quad c(t) :=\sum_{n\ge0)}c_nt^n
\]
then we may rewrite our recurrence as
\[
	t^{-1}(X(t)-x_0) =AX(t) +u(t)b
\]
and consequently
\begin{equation}
	\label{eq:XAb}
	X(t) =(I-tA)^{-1}x_0 +tu(t)(I-tA)^{-1}b.
\end{equation}
Thus we have two distinct contributions to the behavior of the system: one
determined entirely by $A$ and the initial state $x_0$, the other determined by
$A$, $b$ and $u(t)$. It follows from \eqref{eq:XAb} that the state of the system
is always in the column space of the \textsl{controllability matrix}
\[
	W =\pmat{b&Ab&\hdots&A^{d-1}b}
\]
The system is \textsl{controllable} if $W$ is invertible.

(Note that our ``exposition'' of control theory is confined to the simplest case.
In general $b$ and $c$ are replaced by matrices $B$ and $C$. The system is then
controllable if the the $A$-module generated by $\col(B)$ is $\fld^v$, and
observable if the module generated by $\col(C)$ is $\re^v$. This more general
case forced itself on us in our treatment of Laplacians.)

It is convenient to assume $x_0=0$.  Then we have
\[
	c(t) =tu(t)\,c^T(I-tA)^{-1}b.
\]
If the \textsl{observability matrix}
\[
	\pmat{c^T\\ c^TA\\ \vdots\\c^TA^{d-1}}
\]
is invertible, then it is possible to infer the state of the system at time $m$
from the observations $\seq cm{m+d-1}$.  In this case we say that the system is
\textsl{observable}.  Note that the system is observable if and only the pair
$(A,b)$ is controllable.

We can consider a more general version of \eqref{eq:control}: suppose $A$ is
$n\times n$ and $B$ is $n\times k$. We then have a system
\[
	x_{n+1} =A x_n +Bu_n,\quad (n\ge0)
\]
where now $u\in\re^k$. In this case the system is controllable if the $A$-module
generated by the column space of $B$ is $\re^n$. This case arose in
\tref{2dcontrol}.

The series
\[
	c^T(I-tA)^{-1}b,
\]
is known as the \textsl{transfer function} of the system. In control theory our
variable $t$ is normally replaced by a variable $z^{-1}$; thus the transfer
function becomes $c^T(zI-A)^{-1}b$.

%\bibliographystyle{plain}
%\bibliography{/Users/cgodsil/Dropbox/tex/bibs/pst}

\end{document}